\documentstyle[12pt]{article}

\newcommand{\be}{\begin{equation}}
\newcommand{\ee}{\end{equation}}
\newcommand{\bea}{\begin{eqnarray}}
\newcommand{\eea}{\end{eqnarray}}

\newcommand{\ov}{\overline}
\title{Quommutator deformations of spl(N,1)} 
\author{Y. Brihaye \\
Department of Mathematical Physics\\
University of Mons\\
Av. Maistriau, B-7000 Mons, Belgium}
\begin{document}
\begin{titlepage}
\maketitle
\thispagestyle{empty}
\begin{abstract}
We analyse the Witten-Woronowicz's type 
deformations of the Lie superalgebras spl($N,1$)
and obtain deformations parametrized by 
$N(N-1)/2$ independent parameters for $N>2$
and by three parameters for $N=2$. 
For some of these algebras, finite dimensional 
representations can be constructed in terms of
finite difference operators, providing 
operators that are relevant for the classification
of quasi exactly solvable, finite difference equations.
\end{abstract}
\end{titlepage}
\section{Introduction}
The study of deformations of classical Lie algebras
has attracted much attention in the last few years.
Many types of deformations have been proposed. In some of them,
the commutators of the generators are set equal to non linear 
(in fact transcendental)
functions of some of the generators, they are the 
Drinfeld-Jimbo deformations \cite{dri,jim}. Other types of deformations
consist in keeping the right hand sides of the structure relations
linear in the generators while the commutators of the left
hand side are deformed into quommutators, such deformations are called
Witten-Woronowicz type \cite{wit,wor}. 
In this second approach, all the generators are treated on the same footing
and the quommutator relations provide normal ordering rules which allow
to write
the elements of the enveloping algebra in a canonical form. 
This last property is useful for the kind of applications that we have
in mind, namely the study  of the quasi exactly solvable (QES),
finite difference operators. 
\\ The QES (differential or finite difference) 
equations refer to a class of spectral 
equations for which  a finite part of the 
spectrum can be obtained by solving a system
of algebraic equations \cite{tur1,ush}. 
In this respect,
the linear operators which preserve a finite 
dimensional vector space of smooth functions 
are naturally  related to QES equations; such
operators are also called quasi exactly solvable. 
The invariant space is typically a direct sum of 
modules of polynomials of given degrees.  
\\ The observations that QES operators possess hidden
symmetries \cite{tur1,ush} and that they can be related 
to some abstract algebras has motivated several tentatives
to classify these operators and to study their properties
in more details. 
The scalar QES operators of one variable, preserving the 
vector space of polynomials of degree at most $n$, 
are related to the enveloping algebra of the Lie algebra
sl(2) \cite{tur2}. 
Choosing, more generically, $V$ as a direct sum of modules of polynomials,     
it is likely \cite{turshif,gko1,bk,bggk,bn}
that a finite number of basic QES operators 
can be put in correspondance  with the generators of an 
abstract algebra, suitably represented in terms of differential operators.
The more general QES operators can then be obtained
as the polynomials in the basic operators, 
in other words they are the elements of an enveloping algebra. 
\\ The present understanding of QES equations suggests 
that the scalar QES equations are related to
the classical Lie algebras \cite{turshif,gko1}
while the systems of coupled
QES equations are related to superalgebras \cite{turshif,bk,bn}.
For instance the superalgebra osp(2,2) 
(or equivalently spl(2,1)) plays
a crucial role \cite{turshif,bk} for the systems of
two equations of one variable. In the case of systems of two equations
of $V$ variables, the underlying symmetry is related to
the superalgebras spl(V+1,1) \cite{bn} 
\\ On the other hand it seems that, 
if we want to describe the algebraic structure underlying
the QES finite difference equations, the representations
of some deformed algebra emerge in a natural way.
For instance, the scalar, finite difference 
QES operators in one variable 
\cite{tur2} are related to the Witten deformation of sl(2). 
Owing for this example and the results of \cite{bn}  
it is natural to try to classify the deformations of Witten-Woronowicz
type of  the Lie superalgebras spl(N,1).
This is the object of this paper; Drinfeld-Jimbo type 
deformations of general Lie superalgebras are discussed
in \cite{vinet} 
\section{The algebra spl$(N,1)$}
The graded Lie algebra spl$(N,1)$ is generated by
$N^2$ bosonic generators, $N$ fermionic and $N$
anti-fermionic generators. We will denote them
respectively by $E_a^b, V_a$ and $\overline V^a$,
$a,b=1,\dots ,N$. The structure relations read
as follows
\bea 
\label{c1}&\{ V_a , V_b \} &= 0 \ \ \ , \\
\label{c2}&\{ \overline V^a , \overline V^b \} &= 0 \ \ \ , \\
\label{c3}&\{ V_a , \overline V^b \} &= E_a^b  \ \ \ , \\
\label{c4}& [ E_a^b , V_c ] 
&= \delta_c^b V_a - \delta_a^b V_c \ \ \ ,\\
\label{c5}& [ E_a^b , \overline V^c ] &= 
\delta_a^b \overline V^c - \delta_a^c \overline  V^b \ \ \ ,\\
\label{c6}& [ E_a^b , E_c^d ] 
&= \delta_c^b E_a^d - \delta_a^d E_c^b \ \ \ . 
\eea
As indicated in (\ref{c6}), the generators 
$E_a^b$ form a gl($N$) subalgebra. 
Eqs.(\ref{c4}) (resp.(\ref{c5})) show that the generators 
$V_a$'s (resp. $\overline V^b$) form
an $N$-dimensional representation of gl($N$) under
the adjoint action  of the $E_a^b$.  
\\ We will consider the deformations of spl($N,1$)   
obtained by replacing the commutators (\ref{c4})-(\ref{c6}) 
and the anti-commutators (\ref{c1})-(\ref{c3}) 
respectively by quommutators and anti-quommutators~: 
\be
   [A , B]_q = AB - q BA \ \ \ \ , \ \ \ \ \{A , B \}_q = AB + qBA
\label{quomm}
\ee
where $q$ represents a parameter of deformation.
More specifically, we classify the deformations 
obeying the following conditions:
\begin{description}
\item{i)} each commutator (resp. anti-quommutator) 
(\ref{c1})-(\ref{c6}) defining the 
classical algebra is replaced by a quommutator (resp. anti-quommutator)
(\ref{quomm}) with its own parameter $q$.
\item{ii)} the couples of generators which (anti-) commute are imposed to
(anti-~) quommute.
\item{iii)} the right hand side of (\ref{c3})-(\ref{c5})
is kept linear in the generators; in the right hand side of (\ref{c6})
we allow additional terms of the form $V_f \overline V^g$.
\end{description}
\section{Deformations of spl($N,1$)}
In order to present the deformations of spl($N,1$)
it is convenient to introduce a set of  $N(N-1) \over 2$ 
arbitrary parameters labelled $q_{ab}$  
such that
\be
     q_{ba} = { 1 \over q_{ab}} \ \ \ , \ \ \ a,b = 1, \dots , N
\ee
The deformations of spl($N,1$) obeying the conditions above and
compatible with the associativity conditions (i.e. the Jacobi
identities) read as follows~:
\bea
\label{cq1}    &V_a V_b + q_{ab} V_b V_a &= 0 \ \ \ , \\
\label{cq2}    &\overline V^a \overline V^b
               + q_{ab}  \overline V^b \overline V^a &=0 \ \ \ , \\
\label{cq3}   &V_a \overline V^b
               + q_{ba}  \overline V^b V_a &= E_a^b \ \ \ 
\eea
\bea
\label{cq4}   &E_a^b V_c - q_{ac}q_{cb} V_c E_a^b &=
               \delta_c^b V_a - \delta_a^b V_c \ \ \ ,  \\
\label{cq5}   &E_a^b \overline V^c - q_{bc}q_{ca} \overline V^c E_a^b &=
  \delta_a^b \overline V^c - q_{ba} \delta_a^c \overline V_b \ \ \ , 
\eea
\be 
\label{cq6}
      E_a^b E_c^d - q_{ac}q_{cb}q_{bd}q_{da} E_c^d E_a^b =
     \delta_b^c E_a^d - q_{ba}q_{ac}q_{cb} \delta_a^d E_c^b \ \ .
\ee
The normalisation of the generators has been
choosen in such a way that the right hand sides of 
(\ref{cq3}),(\ref{cq4}) look as simple as possible.
The bosonic generators close under the  quommutators (\ref{cq6}),
leading to a deformation of the Lie algebra gl($N$). 
It can be checked that
the products of three and of four $q$'s 
defining the structure constants in (\ref{cq6})
depend effectively on ${(N-1)(N-2) \over 2}$ parameters
among the whole set of the 
${N(N-1) \over 2}$ parameters $q_{ab}$.
In fact, the algebra (\ref{cq6}) coincides with the deformation 
of gl($N$) constructed in \cite{fanu}.
It is amusing to note that the parameters $q_{ab}$,
used to label the deformations in \cite{fanu}, appear
here as defining the anti-quommutators 
(\ref{cq1}),(\ref{cq2}) of the fermionic operators.
Note, however, that another quommutator is used in \cite{fanu},
it can be related to (\ref{cq6}) by an appropriate 
renormalisation of the generators.   
\section{Deformations of spl($2,1$)}
The analysis of Ref.\cite{fanu} reveals that the algebra
gl($2$) leads to a richer set of deformations than for
the generic case gl($N$) for $N>2$. A similar result holds for
 the deformations of spl($N,1$). We now discuss the 
deformations of spl($2,1$) which, remember, is equivalent
to the (perhaps more popular) graded Lie algebra osp($2,2$).
\\ The most general deformation of spl($2,1$) obeying the
three restrictions above
and compatible with associativity depends on three parameters, say
$p,r,s$. The different relations read as follows 
\begin{itemize}
\item for the fermionic-fermionic relations
\be
    V_1 V_1 = V_2 V_2 = 0 \ \ \ , \ \ \ 
   \overline V^1 \overline V^1 = \overline V^2 \overline V^2 = 0
\ee
\be
      \{ V_1,V_2 \}_{{p \over sr}} = 0 \ \ \ , \ \ \ 
      \{ \overline V^1, \overline V^2 \}_{{s \over pr}} = 0 
\label{bff1}
\ee
\be
     \{\overline V^j , V_j \} = E^j_{j} \ \ \ \ , \ \ \ j = 1,2
\label{bff2}
\ee
\be
      \{ \overline V^1, V_2 \}_{psr} = E^1_{2} \ \ \ , \ \ \ 
      \{ \overline V^2, V_1 \}_{{ps \over r}} = E^2_{1} 
\ee
\item for the fermionic-bosonic relations
\bea
&[E^1_{1},V_1] = 0     ,   &[E^2_{2},V_1]_{s^2} = V_1  \nonumber \\
&[E^2_{1},V_1]_{ps/r}=0  ,  &[E^1_{2},V_1]_{sr/p} = -{sr\over p} V_2
 \nonumber  \\
&[E^1_{1},V_2]_{p^2} = V_2  ,    &[E^2_{2},V_2] = 0  \nonumber \\
&[E^2_{1},V_2]_{p/sr}=-{p\over sr} V_1 ,   &[E^1_{2},V_2]_{psr} = 0 
\eea
\bea
&[E^1_{1},\ov V^1] = 0   , 
 &[E^2_{2},\ov V^1]_{1/s^2} = -{1\over s^2} \ov V^1  \nonumber \\
&[E^2_{1},\ov V^1]_{ps/r}= \ov V^2  , 
 &[E^2_{2}, \ov V^1]_{1/psr} = 0  \nonumber \\
&[E^1_{1},\ov V^2]_{1/p^2} = - {1 \over p^2} \ov V^2     , 
 &[E^2_{2},\ov V^2] = 0  \nonumber \\
&[E^2_{1},\ov V^2]_{r/ps}=0  , 
&[E^1_{2}, \ov V^2]_{s/pr} = \ov V^1 
\eea
\item for the bosonic-bosonic operators
\be
    [E^1_{1}, E^2_{2}] = 0 \ \ , \label{b01} \\
\ee
\bea
    &[E^1_{1}, E^2_{1}]_{1/p^2} = - {1 \over p^2} E^2_{1}  ,   
    &[E^2_{2}, E^2_{1}]_{s^2} =  E^2_{1} \label{b0pm} \\   
    &[E^1_{1}, E^1_{2}]_{p^2} =  E^1_{2}  ,   
    &[E^1_{1}, E^1_{2}]_{1/s^2} = - {1 \over s^2}  E^1_{2} 
\label{b1pm}   
\eea
\be
    [E^2_{1}, E^1_{2}]_{s^2/p^2} = E^1_{1} - {s^2 \over p^2} E^2_{2}
 +(s^2 -1)V_1 \overline V^1 - {s^2 \over p^2}(p^2-1) V_2 \overline V^2
\label{bpm}
\ee
\end{itemize}
The three parameters of the deformation, noted $p,r,s$, are 
intrinsic and cannot be eliminated by a rescaling of the generators.
This is seen easily from the way these parameters enter in 
the different coefficients defining the quommutator. 
The undeformed osp(2,2) algebra is recovered by
$p=s=r=1$. 
\\ In the limit $s^2 = p^2 = 1$, the generic case of the previous section
is recovered~: there is only one deforming parameter, noted $r$.
The gl($2$) subalgebra is undeformed.
\\ Keeping $r,s,p$ arbitrary and
decoupling the fermionic generators from  eq.(\ref{bpm}) leads, 
together with (\ref{b01}), (\ref{b0pm}) and (\ref{b1pm}),  to the
two parameters deformation of gl(2) obtained in ref. \cite{faza}
(a four parameters deformation of gl(2) was further 
obtained in ref.\cite{fanu}). 
\section{Representations}
Restricting the three parameters $p,r,s$ in the algebra
presented in sect.~4  to the case 
\be
    p = s^{-1} = q^{{1\over 2}} \ \ \ , \ \ \ r=1
\ee
leads to a one parameter deformation  of osp(2,2) parametrized by $q$
which we will refer to as osp(2,2)$_q$.
It is equivalent, up to a rescaling of the generators,
to the deformation first discussed in \cite{bgk}.
As for the Witten type deformation of sl(2) \cite{tur2},
it is possible to construct some representations 
of osp(2,2)$_q$ in terms
of a finite difference operator $D_q$~:
\be
       D_q f(x) =  { f(x) - f(qx) \over (1-q)x} \ \ \ , \ \ \ 
       D_q x^n = [n]_q x^{n-1} \ \ \ , \ \ \   
       [n]_q \equiv {1 - q^n \over 1-q}  
\label{jack}
\ee 
The representations we want to discuss act on the vector space 
$P(n-1) \oplus P(n)$ (where $P(n)$ denotes the module of
of polynomials of degree  $n$ in the variable $x$).
It is therefore of dimension $2n-1$.
The fermionic generators are represented by
\bea 
     &V_1 = \sigma_-  ,  &V_2 = x \sigma_- \nonumber \\
     &\ov V_1 = q^{-n} (xD_q - [n]_q) \sigma_+  , 
     &\ov V_2 = q^{-1} D_q \sigma_+ 
\label{rep2}
\eea
where $\sigma_{\pm} = (\sigma_1 \pm i \sigma_2)/2$.
The operators $E^i_j$ are easily
constructed with (\ref{bff1}),(\ref{bff2}). 
The  2$\times$2 matrix,  finite difference
operators preserving $P(n-1) \oplus P(n)$ 
are then the elements of the enveloping 
algebra generated by (\ref{rep2}). 
By construction,
these operators are quasi exactly solvable.
\section{Deformations of osp(1,2)}
Recently a Witten type deformation of the super Lie algebra osp(1,2)
was constructed \cite{chung}.
We showed that this deformation is the only one to fulfil
the three requirements of sect. 2, 
its finite dimensional representations can be expressed 
it terms of the operator (\ref{jack}).
Again the space of the representation is 
$P(n-1) \oplus P(n)$ and the two fermionic generators of osp(1,2)
are represented by  2$\times$2 matrix operators.
Using exactly the same notation as in ref.\cite{chung} we find
\be
V_- =
\left(\begin{array}{cc}
0     & D_{q^2}     \\
1     &0  
\end{array}\right)
\ \ \ , \ \ \ \ 
V_+ =
\left(\begin{array}{cc}
0                 &q^{-2n}(x D_{q^2} - [n]_{q^2})     \\
x     &0  
\end{array}\right)
\ee
The three bosonic operators, $H,J_-,J_+$  can then be computed
through the structure of osp(1,2), namely
\be 
    H   =  \{V_-,V_+ \}_q \ \ , \ \ 
    J_- =  \{V_-,V_- \}_q \ \ , \ \
    J_+ =  \{V_+,V_+ \}_q 
\ee
Correspondingly, the Casimir operator
(eq.(13) in ref. \cite{chung})
has a value $ C = -{1\over 2} [-n - {1\over 2}]_{q^2}$
\section{Concluding remarks}
The representations of spl(V+1,1) formulated in terms
of differential operators provide the building blocks for
the construction of the QES operators preserving more
general polynomial spaces of  $V$ variables
\cite{bn}.
This result, and the fact that the quommutators deformations
of super Lie algebras have not been studied systematically
(at least to our knowledge),
motivate a study of the deformations of 
the Lie superalgebras spl(N,1)
by means of quommutators. 
\\ On the other hand,
the most interesting examples of QES systems
(e.g. the relativistic Coulomb problem
\cite{bdk}, the doubly periodic Lame equation
and the stability of the sphaleron in the abelian Higgs model
\cite{bk}) are related to
the algebra osp(2,2).
Generalisations of this algebra,
e.g. its deformations, and the corresponding realizations
in terms of finite difference operators 
therefore desserve some attention.
The operators (\ref{rep2}) could be relevant for the
study of some coupled problems arising in discrete
quantum mechanics in \cite{wiza}.
\\ On a more abstract level we mention that 
 the operators (\ref{rep2}) 
can be used to construct finite difference operators preserving
the vector space $P(m)\oplus  P(n)$. The algebraic structure
underlying these operators is probably 
determined by a series of deformed, non-linear superalgebras
indexed by  $\vert m-n \vert$.
Finally, it would be interesting to relate the two
parameters deformation of osp(2,2) of ref.\cite{para}
with our deformation. The occurence of some mapping between the two
structures would allow to transport the Hopf structure of
\label{para} to our case.

\end{document}